\documentclass[a4paper]{article} %version ecran

\usepackage[english]{babel}

\usepackage[utf8]{inputenc}
\usepackage{amsmath,amsthm}
\usepackage{amsfonts,amssymb}

\usepackage{array}

\usepackage{hyperref}
\usepackage{graphicx}

\usepackage{tikz,tikz-cd} %Graĥique

\usepackage{xifthen} %provide "\ifthenelse" et "\isempty{}"

\usepackage{titlesec}
\titleformat{\subsection}[runin]{\normalfont}{\thesubsection}{0pt}{}[.]

\newcommand{\toposdefaut}{0}
\newcommand{\topos}[1][\toposdefaut]{ 
\ifthenelse{\equal{#1}{0}}{ \mathcal{T} }
{
\ifthenelse{\equal{#1}{1}}{ \mathcal{E} }{ #1 }
}
}

\newcommand{\un}[1][\toposdefaut]{ {1_{\topos[#1]}} }

\newcommand{\Rt}[1][\toposdefaut]{\mathbb{R}_{\topos[#1] }}

\newcommand{\sub}{\textsf{Sub}}
\newcommand{\set}{\text{Set}}

\newcommand{\scal}[2]{ \left\langle #1 , #2 \right\rangle }

\newcommand{\R}{\mathbb{R}}
\newcommand{\C}{\mathbb{C}}

\newcommand{\Ecal}{\mathcal{E}} 
 
\newcommand{\Tcal}{\mathcal{T}}

\newcommand{\Ocal}{\mathcal{O}}

\newcommand{\Hcal}{\mathcal{H}}

\newcommand{\Lcal}{\mathcal{L}}

\newcommand{\Ccal}{\mathcal{C}} 
 
\newcommand{\Bcal}{\mathcal{B}}

\newcommand{\block}[1]
{

\par \subsection{} #1

\bigskip}

\newcommand{\blockn}[1]{\par #1 \bigskip}

\newcommand{\Th}[1]
	{
	\bigskip	
	\textbf{Theorem : }{\itshape #1}
		
	\bigskip
	}

\newcommand{\Prop}[1]
	{

	\bigskip
	
	\textbf{Proposition : }{\itshape #1}
		
	\bigskip
	
	}

\newcommand{\Cor}[1]
	{

	\bigskip
	
	\textbf{Corollary : }{\itshape #1}	
		
	\bigskip

	}

\newcommand{\Lem}[1]
	{

	\bigskip
	
	\textbf{Lemma : }{\itshape #1}
		
	\bigskip
	
	}

\newcommand{\Def}[1]
	{
	
	\bigskip
	
	\textbf{Definition : }{\itshape #1}
	
	\bigskip
	
	}

\newcommand{\Dem}[1]{
	
	\smallskip
	
	\textbf{Proof : } \par
	 {#1} $\square$
	 
	 \bigskip
}

\begin{document}

\selectlanguage{english}
\pagestyle{plain}
\title{Measure theory over boolean toposes}
\author{Simon Henry}
%\address{Simon.~Henry: \\ Paris, F-75005 France}
%\email{henry\@@phare.normalesup.org}
%\keywords{mono\"{\i}ds, hypergroups}
%\subjclass[2010]{20N20,43A62} 

\renewcommand{\thefootnote}{\fnsymbol{footnote}} 
\renewcommand{\thefootnote}{\arabic{footnote}} 

% 18B25 : topoi
% 03G30 : Categorical logic, topoi
% 46L10 : General theory of von Neumann algebras
% 46L51 : Noncommutative measure and integration

\maketitle

\begin{abstract}
In this paper we develop a notion of measure theory over boolean toposes which is analogous to noncommutative measure theory, i.e. to the theory of von Neumann algebras. This is part of a larger project to study relations between topos theory and noncommutative geometry. The main result is a topos theoretic version of the modular time evolution of von Neumann algebra which take the form of a canonical $\mathbb{R}^{>0}$-principal bundle over any integrable locally separated boolean topos.

\end{abstract}

\tableofcontents

\section{Introduction}

\blockn{This paper is part of a larger project to understands the relation between noncommutative geometry and topos theory, and is more precisely focused on the measure theoretic aspect of this relation.}

\blockn{An extremely efficient way to relate a topos $\Tcal$ to objects from noncommutative geometry is the study of internal Hilbert spaces of $\Tcal$. Let us explain a little more this point. In the same way that an abelian category is a category which behave like the category of abelian groups, a topos is a category that behave essentially like the category of sets. In the case of abelian category, it is known that any result that can be proved in the category of abelian groups by diagram chasing will hold in any abelian category. In the same way, it is know that any result that can be proved in the category of sets without using neither the law of excluded middle nor the axiom of choice will holds in any topos. This interpretation is essentially what is called the internal logic of a topos, and it allows, for any topos $\Tcal$, to transport, for example, the definition of an Hilbert space into a definition of an internal Hilbert space in $\Tcal$. An Hilbert space of $\topos$ (or a $\topos$-Hilbert space) is an object of $\topos$ endowed with a series of operation satisfying the axiom for an Hilbert spaces in the internal logic of $\topos$. It appears that\footnote{We might need $\Tcal$ to be a Grothendieck topos and not just an elementary topos for this.} the category of $\topos$-Hilbert spaces and bounded operators between them is an ``external" $C^{*}$-category (in the sense of for example \cite[Definition 1.1]{wstarcat}) naturally attached to the topos $\Tcal$ and related to its geometry.}

\blockn{When the internal logic of a topos $\Tcal$ satisfies the law of excluded middle ($\Tcal$ is then said to be a boolean topos) then one can, using internal logic, construct the supremum of a bounded directed net of positive operators on an internal Hilbert space simply by constructing their supremum internally (the law of excluded middle is all we need to construct supremum of such families, and one easily see that the internal supremum will also be an external supremum), this turn the $C^{*}$-category of $\Tcal$-Hilbert space into a ``Monotone complete $C^{*}$-category", in the sense that all the $C^{*}$-algebras of endomorphisms are monotone complete $C^{*}$-algebras (see \cite[Definition III.3.13]{takesaki2003theoryI}).}

\blockn{We will show that this monotone complete $C^{*}$-category is a $W^{*}$-category if and only if $\Tcal$ satisfies a condition of existence of measures called ``integrability" (see definition \ref{DEfintegrable} and theorem \ref{THintegrable}).}

\blockn{For these reasons, it seem natural to think of boolean toposes as the topos theoretic analogue of monotone complete $C^{*}$-algebras and of integrable boolean topos as the analogue of von Neumann algebras. This is encouraged by the fact that the ``commutative case" agrees: the category of boolean locales is well known to be anti-equivalent to the category of commutative monotone complete $C^{*}$-algebras and normal morphisms between them; and, under this correspondence, commutative von Neumann algebra correspond exactly to the integrable boolean locale.

This being said we acknowledge the fact that there do exist non-boolean topos whose category of internal Hilbert space is a $W^{*}$-category, but it seems to us that these are examples of toposes whose geometry is not reflected by their $C^{*}$-category of internal Hilbert spaces, and that they should be discarded of the present work.}

\blockn{For these reasons, and as we are concerned with measure theory, we will only focus on boolean toposes. \emph{All the topos mentioned in this paper are boolean Grothendieck toposes.}}

\blockn{Von Neumann algebras (and also $W^*$-categories) are well known\footnote{see for example \cite[Chapter VIII]{takesaki2003theoryII} for von Neumann algebras and \cite[section 3]{wstarcat} for the case of $W^{*}$-categories.} to have a canonical ``modular" time evolution. The main result of this paper is to provide a geometric description of this time evolution in the case of von Neumann algebras arising from a boolean locally separated topos $\topos$ in terms of a certain (canonical) principal $\Rt^{>0}$-bundle over $\topos$ analogous to the bundle of positive locally finite well supported measures over a boolean locale (which is a principal bundle exactly because of the Radon-Nikodym theorem).}

\blockn{
This gives a classification in types I, II and III of boolean integrable locally separated toposes analogous (but not totally equivalent) to the classification in type of von Neumann algebras: Type I corresponds to separated toposes, type II to toposes which are not separated but which have a trivial modular bundle, and type III to toposes which have a non-trivial modular bundle.

Of course in full generality, one obtains that every boolean integrable locally separated topos decomposes in a disjoint sum of a topos of each of the three types by applying this disjunction internally in its localic reflection.

}

\blockn{In the last section of this paper, we consider $\topos$ a boolean integrable locally separated toposes, $X$ an object of $\topos$ such that $\topos_{/X}$ is separated and $l^2(X)$ the Hilbert space of square summable sequences. In this situation, we show that the modular time evolution of the von Neumann algebra $A$ of globally bounded endomorphisms of $l^2(X)$ is indeed the time evolution on $A$ described by $\chi$, and that an invariant measure on $\topos$ induces a trace on $A$. }

\blockn{The following table sum-up the dictionary between topos theory (in the left column) and operator algebra (in the right column) that arise in this paper. This is of course just a vague analogy that we have observed while developing this theory, we do not claim that there is any sort of rigorous correspondence here. In particular we think that this dictionary is meant to be made more precise in the future. For example a work in progress (mentioned in the introduction of the author's thesis under the name ``Non abelian monoidal Gelfand duality") highly suggest that it can be made a lot more precise if we take into account on the right hand side the monoidal structures that arise from the internal tensor product of Hilbert space and the compatibility to these structures. }

\renewcommand{\arraystretch}{1.5}
\blockn{

\begin{center}
\begin{tabular}{|>{\centering}m{5.5cm} |m{5.5cm}<{\centering}|}

\hline 
 
Boolean topos (locally separated) & Monotone complete $C^{*}$-algebra (up to morita equivalence) \\

\hline

Boolean integrable topos (locally separated) & $W^{*}$-algebra (up to morita equivalence) \\

\hline

Localic reflection & Center \\

\hline

Measure on an object & Semi-finite normal weight \\

\hline

Measure on an object of mass $1$ & Normal state \\

\hline

Invariant measure & Normal (semi-finite) trace \\

\hline

The modular bundle & The $\Delta$ operator of the Tomita-Takesaki construction \\

\hline

The family of line bundles $(F_t)_{t \in \R}$ & The modular time evolution \\

\hline

\end{tabular}
\end{center}
}
\renewcommand{\arraystretch}{1.0}

\blockn{This paper essentially corresponds to the second chapter of the author thesis (\cite{mythesis}). Although the presentation of the present paper has been revised, and the proof technique improved in comparison to the version present in the author's thesis.}

\section{Preliminaries and notations}

\block{The letter $\topos$ will denotes a boolean (Grothendieck) topos. Its terminal object is denoted $\un$, if $X$ is an object of $\topos$ we denote by $\topos_{/X}$ the slice topos whose objects are couples $(Y,f)$ where $f$ is a morphism from $Y$ to $X$. It corresponds geometrically to the etale space of $X$. A family $(X_i)$ of object of a topos is called a generating family or a family of generator if any object of $\topos$ can be covered by epimorphic image of the $X_i$ (in which case on can construct a site for $\topos$ on the full sub-category of $\topos$ whose objects are the $X_i$).}

\block{$\Rt$ denote the object of real number of the topos $\topos$. As $\topos$ is boolean, all reasonable\footnote{except the definition by Cauchy sequences which require the axiom of depend choice to be equivalent to the others.} definitions of the real numbers agree and give rise to an internal ordered fields satisfying the property that every bounded set has a supremum.
}

\block{If $X$ is an object of $\topos$ we denote by $|X|$ the internal cardinal of $X$, which is a function from the terminal object of $\topos$ to $\mathbb{N} \cup \infty$ (i.e. a function to the object $\mathbb{N} \cup \infty$ of $\topos$ or equivalently a map from the locale associated to the frame of sub-object of $\un$ to the discrete locale whose points are $\mathbb{N}\cup \infty$). }

\block{An object of a topos will be called \emph{finite} if it is internally finite (either Kuratowski finite or cardinal finite which are equivalent in boolean topos). A map $f:X\rightarrow Y$ between object of a topos is said to be finite if internally one has $\forall y \in Y, f^{-1}(\{y\})$ is finite. It is said to be $n$-to-$1$ if internally one has $\forall y \in Y, f^{-1}(\{y\})$ has cardinal $n$, one can eventually take $n = \infty$ in this definition.}

\block{\label{valuation}Let $B$ be a complete boolean algebra with $0$ and $1$ the bottom and top elemnts. A valuation on $B$ is a non-decreasing function $\mu$ from $B$ to the set $\R_+^{\infty}$ of possibly infinite non-negative real numbers, such that $\mu(U)+\mu(V)= \mu(U \cup V) + \mu(U \cap V)$, $\mu(0)=0$ and $\mu$ preserve arbitrary directed supremum.

A valuation is said to be finite if $\mu(1)<\infty$, locally finite if $1$ is a union of element of finite valuation and well supported if the only element of valuation $0$ is $0$.

}

\block{\label{integral}If $X$ is a boolean locale a measure on $X$ is valuation on the complete boolean algebra $\Ocal(X)$. If $f : X\rightarrow \C$ is a continuous map and $\mu$ a valuation on $X$ one can integrate $f$ against $\mu$. (see for example \cite{vickers2008localic} for a general treatment of the integration over locale in constructive mathematics). This integral satisfy all the properties of the usual notion of integral: it is linear, additive and preserve directed supremum when applied to positive functions. For general function, we need to restrict to function for which the integral of $|f|$ is finite and then it is linear and satisfy the various form of Lebesgue's dominated convergence theorem\footnote{This last fact can be deduced from the constructivity of Vickers' theory: it allows to deduce directly lower semi-continuity of integrals with parameters, and using to the domination hypothesis one obtains the lower semi-continuity of $-\int$ and hence the continuity of the integral.}.}

\block{Boolean locales satisfy the Radon-Nikodym theorem: If $\mu$ and $\nu$ are two locally finite well supported measures on a boolean locale $X$, then there exists a (unique) continuous map $f$ from $X$ to $\R^{>0}$ such that $\mu= f.\nu$. This last equality mean that for all open sublocales $U \subset X$ one has $\mu(U) = \int f \mathbb{I}_U d\nu$, see for example \cite[section 3.4]{jackson2006sheaf} for more details.}

\section{Hilbert spaces and integrable toposes}
\label{Integrabletopos}

\block{As mentioned in the introduction, we call an internal Hilbert space or a $\Tcal$-Hilbert space an object of $\topos$ endowed with operations making it into an internal Hilbert space. A morphism of $\Tcal$-Hilbert space, is a map $f$ which is internally a bounded operator, and such that the norm is externally bounded, i.e. such that there exists an external constant $K$ such that internally $\Vert f \Vert \leqslant K$ (we also call them ``globally bounded operator"). We denote $\Vert h \Vert_{\infty}$ the smaller such constant $K$, while $\Vert h \Vert$ will denote the internal norm (which is a function on the terminal object). 

Endowed with this norm and the internal adjunction, the category of $\Tcal$-Hilbert space is a $C^{*}$-category\footnote{See for example \cite[Definition 1.1]{wstarcat} for the definition of $C^{*}$-category.}. Moreover, as $\Tcal$ is boolean, it is monotone complete as mentioned in the introduction.
}

\block{If $X$ is an object of a topos, then $l^{2}(X)$ denote the $\Tcal$-Hilbert space a square sumable sequences of complex numbers indexed by $X$. Equivalently, $l^{2}(X)$ is the Hilbert space generated by elements $e_x$ for each $x \in X$ such that $\scal{e_x}{e_{x'}} = 1$ if $x=x'$ and $0$ otherwise. }

\block{Let $\topos$ be a boolean topos, and $X$ be an object of $\topos$, and $\sub(X)$ the complete boolean algebra of sub-object of $X$. We call \emph{measure} on $X$ a valuation on the complete boolean algebra $\sub(X)$ as in \ref{valuation}.

If $f:X\rightarrow \R^+_{\topos}$ is a morphism from $X$ to the object of positive real numbers of $\topos$, then $f$ corresponds externally to a continuous function from the locale associated to the frame $\sub(X)$ to the locale of real numbers. In particular, if $\mu$ is a measure on $X$ one can define:

\[\int_{X} f d\mu \hspace{20pt} \text{  Or  } \hspace{20pt} \int_{x \in X} f(x) d\mu \]

as the integral of this corresponding function with respect to $\mu$ (as in \ref{integral}).
}

\block{\label{lemma_statefrommeasure}This external integration of internal functions allows to relates measures on objects to states of the $C^{*}$ category of the topos:

\Lem{Let $X \in \Tcal$ be an object of $\Tcal$, let $\Hcal$ be an internal Hilbert space of $\Tcal$, $\mu$ a measure on $X$ and $v:X \rightarrow \Hcal$ a function in $\Tcal$ such that:

\[ \int_{x \in X} \Vert v(x) \Vert^{2} d\mu = 1 \]

then, for $f: \Hcal \rightarrow \Hcal$ a globally bounded operator, the formula:

\[\eta(h) := \int_{x \in X} \scal{v(x)}{h v(x)} d\mu \]

defines a normal state $\eta$ on the monotone complete $C^{*}$-algebra $B(\Hcal)$ of globally bounded operators on $\Hcal$.
}

\Dem{ $\eta$ is well defined and of norm smaller than one because:

\[ \int_{x \in X} |\scal{v(x)}{h v(x)}| d\mu \leqslant \Vert h \Vert_{\infty} \int_{x \in X} \Vert v(x) \Vert = \Vert h \Vert_{\infty} \]

the linearity and the positivity are immediate by linearity and positivity of the scalar product and the integral. The value of of $\eta(1)$ follows from the condition relating $v$ and $\mu$. 

The normality of $\eta$ follow from the fact that the supremums are computed internally and that the integral defined in \cite{vickers2008localic} commute to directed supremums.

}

}

\block{\label{DEfintegrable}\Def{A boolean topos $\topos$ is said to be \emph{integrable} if for all non-zero object $X \in \topos$ there exists a non zero finite measure on $X$.}}

\block{\label{THintegrable}\Th{Let $\topos$ be boolean topos, then the $C^{*}$-category of $\topos$-Hilbert spaces is a $W^{*}$-category if and only if $\topos$ is integrable.}

\Dem{Assume first that $\topos$ is integrable. As the category of $\topos$-Hilbert spaces has bi-product, it is enough to show that the algebra of endomorphism of any internal Hilbert space is a $W^{*}$-algebra (because of \cite[Proposition 2.6]{wstarcat}), and as it is already know that these are monotone complete, we only have to prove that they admit enough normal state (see \cite[Theorem III.3.16]{takesaki2003theoryI}).

Let $h : \Hcal \rightarrow \Hcal$ be a non zero positive self adjoint operator, let $X \subset \Hcal$ the set of $x$ such that $\scal{x}{h(x)}>0$ and $\Vert x \Vert \leqslant 1$. If $X=0$ then $h=0$ hence, $X$ is a non zero subset of $\Hcal$ which hence admit a positive finite measure. After rescaling, lemma \ref{lemma_statefrommeasure} gives a state $\eta$ such that $\eta(h)>0$ which concludes the proof of the first implication.

Conversely, assume that the category of Hilbert spaces is a $W^{*}$-category, let $X \in \topos$ a non-zero object of $\topos$, then there is a state $\eta$ on the algebra of endomorphisms of $l^{2}(X)$. For any subobject $S \subset X$ one can define $P_S$ to be the endomorphism of $l^{2}(X)$ defined by $P_S(e_x)=e_x$ if $x \in S$ and $0$ otherwise. Then $\mu(S) = \eta(P_S)$ defines a non zero finite measure on $X$. 
}

}

\block{ \Cor{Let $\topos$ be an integrable boolean topos, then there exists a von Neumann algebra $A$ (uniquely determined up to Morita equivalence of von Neumann algebras) such that the category of $\topos$-Hilbert spaces is equivalent to the category of self-dual Hermitian $A$-modules (in the sense of \cite[example 1.4]{wstarcat}).}

\Dem{As the category of $\topos$-Hilbert spaces has arbitrary orthogonal sums and splitting of projections (they are computed internally) it is enough, by  proposition $7.6$ of \cite{wstarcat}, to check that it has a generator in the sense of \cite{wstarcat} $7.3$.

We fix $G$ a set of generators of $\topos$.
Let $E$ be the family of all isomorphisms class of triples $(H,g,f)$ where $H$ is a $\topos$-Hilbert space, $g$ is an element of $G$ and $f$ is a map $g \rightarrow H$ such that internally the image of this map spans a dense subspace of $H$.
To an element $(\Hcal,g,f)$ of $E$ one can associate a continuous function on $g \times g$ defined by $(x,y) \rightarrow \scal{f(x)}{f(y)}$, and if $(H,g,f)$ and $(H',g,f')$ define the same function on $g \times g$ then they are isomorphic. In particular the isomorphism class of elements of $E$ form a set.

Now, for any $f :A \rightarrow B$ an operator between two $\topos$-Hilbert spaces, there exists $g \in G$ and a map $\lambda : g \rightarrow A$ such that $f \circ \lambda \neq 0$, the adherence of the span of $\lambda$ gives an element $H$ of $E$ and a map $i$ from $H$ to $A$ such that $f \circ i \neq 0$. This proves that elements of $E$ form a family of generators of $\Hcal(\topos)$, i.e. their orthogonal sum is a generator of this $W^*$-category. This concludes the proof.
}

Even if this von Neumann algebra is naturally attached to the topos $\topos$ and uniquely determined by it up to Morita equivalence, it is in most case ``too big" for example, if $\topos$ is the topos of $G$-sets for some discrete group $G$, then $\Hcal(\topos)$ is the category of unitary representations of $G$ and the associated von Neumann algebra is the enveloping von Neumann algebra of the maximal $C^*$-algebra of the group, which is an enormous algebra.

It seems that a more reasonable algebra to consider in practice is the algebra of operators on a space $l^2(X)$ for $X$ a separating bound of $\topos$. In the case of $G$-sets, this gives the Von Neumann algebra of the group, eventually up to Morita equivalence depending on the choice of the object $X$. The results of the last section suggest this algebra can be controlled by the geometry of $\topos$, whereas in general an algebra of operators on an arbitrary Hilbert space over a topos $\topos$ can have nothing in common with the geometry of $\topos$. For example, any von Neumann algebra arises as the algebra of globally bounded endomorphisms of some representation of a discrete group.

}

\section{Invariant measures}

\block{\label{invmeasDef}\Def{Let $\topos$ be a boolean topos. An invariant measure\footnote{It should probably be more correct to call this a ``well supported locally finite invariant measure".} on $\topos$ is a function which to every object $X$ of $\topos$ associates a real number $\mu(X) \in [0,\infty]$ such that:

\begin{enumerate}
\item[(IM1)] For each $X \in \topos$, the restriction of $\mu$ to sub-objects of $X$ defines a locally finite valuation on $\sub(X)$.

\item[(IM2)] There exists a generating family of $\topos$ of objects $X$ on which the valuation induces on $\sub(X)$ is well supported.

\item[(IM3)]If $f:Y\rightarrow X$ is a $n$-to-$1$ map in $\topos$ and if $\mu(Y)<\infty$ then:

\[ \mu(X) = \frac{\mu(Y)}{n} \]

\end{enumerate}
}

Note that the third condition implies that $\mu(X)= \mu(Y)$ when $X$ and $Y$ are isomorphic. In the case $n=\infty$, the axiom $(IM3)$ has to be interpreted as $\mu(X)=0$ if $\mu(Y)<\infty$.
}

\block{This definition, which might seems ad-hoc, essentially come from the study of explicit example of topos where the associated Von Neumann algebra have an explicite time evolution and found is justification in the results of section \ref{sectionmain}. Some detail on its origin can be found in the author's thesis \cite{mythesis} in chapters one and two. }

\block{\label{IMintegral}The following proposition largely extend the strength of axiom $(IM3)$.

\Prop{Let $\topos$ be a boolean topos endowed with an invariant measure $\mu$. Let $f:Y\rightarrow X$ be a finite map and $h$ be a complex valued function on $Y$ (i.e. a function from $Y$ to the object of complex number of $\topos$).

Assume that $h$ is positive or that: 
\[\int_{y \in Y} |h(y)|d\mu < \infty \]

Then one has:

\[ \int_{y \in Y} h(y) d\mu = \int_{x \in X} \left( \sum_{y \in f^{-1}(y)} h(y) \right) d\mu \]

}

\Dem{We will proceed in a series of steps:

\begin{enumerate}

\item Assume first that $h$ is the constant equal to one function, that $\mu(Y)<\infty$ and that $f$ is a $n$-to-$1$ map for some (external) integer $n$. Then the results is exactly $(IM3)$ applied to $f:Y \rightarrow X$.

\item Still assume that $h=1$ and $\mu(Y)<\infty$, but $f$ is an arbitrary finite map. Let:
\[ X_n=\{x \in X| \text{$x$ has exactly $n$ antecedent by $f$  } \}\]
\[ Y_n = f^{-1}(X_n),\]

and let $f_n$ be the restriction of $f$ as a map $Y_n$ to $X_n$. One has $X=\coprod X_n$ and $Y=\coprod Y_n$. For all $n>0$ one can apply the first point to the function $f_n$ to get that $\mu(Y_n)=n\mu(X_n)$. This formula also holds for $n=0$ as $Y_n$ is empty by definition. Using the decomposition of $X$ and $Y$ into the $X_n$ and $Y_n$ one obtain that the left hand side of the formula is $\sum \mu(Y_n)$ and the right hand side is $\sum n \mu(X_n)$ which proves they are equals.

\item Assume now that $Y$ is arbitrary and that $h$ is the characteristic function of $S \subset Y$ with $\mu(Y)<\infty$. Then the function $\sum_{f(y)=x} h(y)$ is supported on $f(S)$ hence the result follows from the previous point applied to $f:S\rightarrow f(S)$.

\item By linearity, this automatically imply the result for any linear combination of such characteristic functions. And as every positive function is a directed supremum of such linear combination and that both the integral and the sum commute to directed supremum of positive function it imply that the result holds for any positive function $h$.

\item if $h$ is an arbitrary complex function such that $\int |h|<\infty$ then we can decompose it into $a-b+ic-id$ where $a,b,c,d$ are four positive functions, and the result follows by linearity of the integral and of the sum.

\end{enumerate}

}

}

\block{\label{IMintegral2}\Cor{Let $\topos$ be a topos endowed with an invariant measure and let $f:Y \twoheadrightarrow X$ be an arbitrary epimorphism, with $\mu(Y)<\infty$ then:

\[ \mu(X) = \int_{y\in Y} \frac{1}{|f^{-1}(f(y))|} d \mu \]

}

\Dem{If $f$ is a finite map then the result follow directly from \ref{IMintegral} applied to $h=\frac{1}{|f^{-1}(f(y))|}$. In the general case, we decompose $X$ and $Y$ into two components: one where $f$ is finite and the other where $f$ is $\infty$-to-$1$. on the first component the previous result apply, one the second component $\mu(X)=0$ because of $(IM3)$ and the integral is also $0$ because $h=0$. }
}

\block{\label{IMextension}\Prop{Let $\Ccal$ be a class of object of $\topos$ such that if $C \in \Ccal$ and $f:X\rightarrow C$ is any map in $\topos$ then $X \in \Ccal$. Also assume that every object of $\topos$ can be covered by objects in $\Ccal$.

Let $\mu$ be a function which associate to every object in $\Ccal$ a non-negative possibly infinite real number, and which satisfy all the axiom of the definition of an invariant measure when restricted to object and map in $\Ccal$. Then $\mu$ extend to a unique invariant measure on $\topos$.
}

We insist on the fact that any subobject of an object in $\Ccal$ is also in $\Ccal$, hence axiom (IM1) do mean something when we restrict to $\Ccal$.

\Dem{First observe that the results of \ref{IMintegral} and \ref{IMintegral2} still holds for our function $\mu$ when $X,Y\in \Ccal$. Let $\Ccal_f$ the subclass of $C \in \Ccal$ such that $\mu(C)<\infty$.

Let $X$ be an object of $\topos$ and assume that there is two epimorphisms $f:C\twoheadrightarrow X$ and $f':C' \twoheadrightarrow X$ with $C,C' \in \Ccal_f$. Then one has:

\[ \int_{c \in C} \frac{1}{|f^{-1}f(c)|}d\mu = \int_{c \in C'} \frac{1}{|f'^{-1}f'(c)|}d\mu   \] 

Indeed, to see that, let $P=C \times_X C'$ and $g$ and $g'$ the two maps from $P$ to $C$ and $C'$, let $h$ be the function on $P$ defined by $h(c,c')= \frac{1}{|f'^{-1}f'(c')|} * \frac{1}{|f^{-1}f(c)|}$. Now, (internally) for any $c \in C$ one has:

\[ \sum_{x \in P \atop g(x)=c} h(x) = \sum_{c' \in C' \atop f'(c')=f(c)} \frac{1}{|f'^{-1}f'(c')|} \frac{1}{|f^{-1}f(c)|} = \frac{1}{|f^{-1}f(c)|}\]

And the same computation exchanging $C$ and $C'$ together with proposition \ref{IMintegral}, yields that:

\[ \int_{c \in C} \frac{1}{|f^{-1}f(c)|}d\mu = \int_P h =  \int_{c' \in C'} \frac{1}{|f'^{-1}f'(c')|}d\mu \]

In particular, for any object $X$ of topos which is a surjective image of $C \in \Ccal_f$ there is a uniquely defined extension of $\mu$ to $X$ as:

\[ \mu(X)= \int_{c \in C} \frac{1}{|f^{-1}f(c)|}d\mu \]

If $S \subset X$ then $S$ is the surjective image of $f^{-1}(S)$ hence one also have a definition of $\mu(S)$ and this defines a valuation on $X$. If $X$ and $X'$ are two subobjects of $Y$ begin epimorphic images of objects in $\Ccal_f$, then the measure defined on them agree on their intersection by the uniqueness property, hence this measure extend as a measure on their union, and as every object is covered by surjective image of object of $\Ccal_f$ by hypothesis this defines a unique measure on each object of $\topos$. We already proved that the extension satisfy (IM1), and it clearly satisfy (IM2) because $\Ccal$ already contains a generating family of object on which the measure is well supported.

In order to conclude, we just need to prove that this extension satisfy (IM3). Let $f:Y \rightarrow X$ be a $n$-to-$1$ map. Let $C\in \Ccal_f$ and $v:C\rightarrow X$. Let $C':= C\times_X Y$. The map $\pi_1: C' \rightarrow C$ is again $n$-to-$1$ hence: 

\[ \mu(v(C)) = \int_{c \in C} \frac{1}{|v^{-1}v(c)|} d\mu = \frac{1}{n} \int_{c' \in C'} \frac{1}{|v^{-1}v(\pi_1(c'))|} d\mu. \]

On the other hand:

\[ \mu(\pi_2(C'))= \int_{c' \in C'} \frac{1}{|\pi_2^{-1}(\pi_2(c')|} d\mu \]

But if $c'=(c,y)$ then $\pi_2^{-1}(\pi_2(c'))$ is the set of $(a,y)$ with $a\in v^{-1}v(c)$ hence, $|\pi_2^{-1}(\pi_2(c'))|=|v^{-1}v(\pi_1(c'))|$ so that:

\[ \mu(v(C))= \frac{1}{n} \mu(\pi_2(C')) =   \frac{1}{n} \mu(f^{-1}(v(c)))\]

As the measure on $X$ is defined out of the measure on such $v(C)$ this shows that $\mu(X)=\frac{1}{n} \mu(Y)$ and concludes the proof.

}

}

\block{Let $\topos$ be a topos, if $X$ is an object of $\topos$ we denote by $M(X)$ the set of invariant measures on the slice topos $\topos_{/X}$. If $f:Y\rightarrow X$ is a map in topos, and $\mu \in M(X)$ then we define $f^{*}\mu \in M(Y)$ by the formula $f^{*}(\mu)(V) = \mu(f_! V)$. Where $f_!$ is the ``composition with $f$" functor from $\topos_{/Y}$ to $\topos_{/X}$. This defines a contravariant functor from $\topos$ to $\set$.

\Th{Let $\topos$ be a boolean topos, then the contravairant functor $M$ of invariant measure is representable by an object denoted $\chi$.
}

\Dem{As $\topos$ is a Grothendieck topos we just need to prove that $M$ is a sheaf for the canonical topology of $\topos$. More precisely: let $f: X \twoheadrightarrow Y$ be an epimorphism, $P=X \times_Y X$ and $\pi_1,\pi_2: P \rightrightarrows X$ the two projections, we need to prove that for any $\mu \in M(X)$ such that $\pi_1^{*} \mu = \pi_2^{*}\mu$ there exists a $\nu \in M(Y)$ such that $\mu=f^{*}\nu$.

By working with $\topos_{/Y}$ instead of $\topos$ one can freely assume that $Y$ is the terminal object of $\topos$ and $P=X \times X$.

Let $\Ccal$ be the class of object $V$ of $\topos$ such that there exists a map from $V$ to $X$. Let $V \in \Ccal$ and let $f,g:V \rightrightarrows X$ two maps from $V$ to $X$, then $(V,f,g)$ is an object of $\topos_{/P}$ one has by hypothesis $\pi_1^{*}\mu (V,f,g) = \pi_2^{*}\mu(V,f,g)$, but one the other hand, $\pi_1^{*}\mu (V,f,g) = \mu(V,f)$ and $\pi_2^{*}\mu(V,f,g) = \mu(V,g)$. Hence $\mu(V,f)$ does not depends on $f:V\rightarrow X$ and we define, for $V \in \Ccal$,  $\nu(V)= \mu(V,f)$ for any $f:V\rightarrow X$. It is easy to check that this satisfy the axiom of an invariant measure restricted to $\Ccal$: for each axiom, one can chose a coherent family of function to $X$ and apply the corresponding axiom for the measure $\mu$ in $\topos_{/X}$. Hence by proposition \ref{IMextension} this extend to a unique measure $\nu$ on $\topos$ which by definition satisfy $f^{*}\nu= \mu$
}

}

\block{Of course, this object $\chi$ comes with several structures (due to structures on $M$). The most important is the multiplicative action of $\Rt^{>0}$ on $\chi$ which corresponds to the fact that if $\mu$ is an invariant measure on $\topos_{/X}$ and $f$ a function from $X$ to $\Rt^{>0}$ then one can define, exactly as in the case of classical measure theory, $f.\mu$ by $f.\mu(Y,p) = \int f(p(y)) d\mu$. One can also define an addition on $\chi$, but there is no ``zero" element because it is ruled out by (IM2).}

\section{Separated and locally separated boolean toposes}

\blockn{In this section we introduce the notion of separated and locally separated toposes due to I.Moerdjik and C.J.Vermeulen in \cite{moerdijk2000proper}, and we prove a new characterization separation for boolean toposes (see Theorem \ref{sep=locfinit}). }

\block{The following definitions are due to Moerdijk and Vermeulen (see \cite{moerdijk2000proper})
\Def{
\begin{itemize}
\item A topos is said to be compact if its localic reflection is compact.
\item A geometric morphism $f:\topos[1] \rightarrow \topos$ is said to be proper if internally in $\topos$, the $\topos$-topos $\topos[1]$ is compact.
\item A geometric morphism $ f: \topos[1] \rightarrow \topos$ is said to be separated if its diagonal map $\Delta: \topos[1] \rightarrow \topos[1] \times_{\topos} \topos[1]$ is a proper map.
\item A topos is said to be separated if its canonical geometric morphism to the base topos is separated.
\end{itemize}
}

Being separated is a strong property, for example the topos of $G$-sets for a discrete group $G$ is separated if and only if $G$ is finite.
}

\block{\label{sep=locfinit}The main result of this section is the following characterization of separation in the case of boolean toposes:

\Th{Let $\topos$ be a boolean topos. Then $\topos$ is separated if and only if $\topos$ has a generating family of internally finite object.}

The proof of this theorem will be dived in various lemma, one implication is proved in corollary \ref{sepLEM1} and the other by corollary \ref{sepLEM2}.

}

\block{\label{hyperconexseparated}\Prop{Let $\topos$ be a hyperconnected\footnote{This means that its localic reflection is the point, i.e. that $\un$ has no non-trivial subobjects, see \cite[A4.6]{sketches}.} separated topos. Then $\topos$ is atomic and all its atoms are (internally) finite.}

\Dem{Let $\topos$ be a hyperconnected topos. Let $\Bcal$ be a non zero boolean locale endowed with a geometric morphism to $\topos$. In the logic of $\Bcal$, the pullback of the topos $\topos$ is hyperconnected separated and has a point. Applied internally in $\Bcal$, the theorem II.3.1 of \cite{moerdijk2000proper} shows that the pullback of $\topos$ is equivalent to the topos of $G$-sets for $G$ a compact localic groups in $\Bcal$. In particular the pullback of $\topos$ is atomic in the logic of $\Bcal$ and as the map from $\Bcal$ to the base topos is an open surjection (because the base topos is boolean) it implies that $\topos$ is atomic in the base topos (by \cite[C5.1.7]{sketches}).

Let now $a$ be an atom of a separated topos $\topos$, then $\topos_{/a} \rightarrow *$ is a proper map (because as $a$ is an atom it is hyperconnected), and $\topos \rightarrow *$ is separated by hypothesis. Hence, by Proposition II.2.1(iv) of \cite{moerdijk2000proper} this proves that the map $\topos_{/a} \rightarrow \topos$ is proper, i.e. that $a$ is internally finite.
}

}

\block{\label{sepLEM1}We are now ready to prove the first half of theorem \ref{sep=locfinit}.

\Cor{Let $\topos$ be a boolean separated topos, then $\topos$ admit a generating family of internally finite objects.}

\Dem{Let $\Lcal$ be the localic reflection of $\topos$. The geometric morphism from $\topos$ to $\Lcal$ is hyperconnected by definition and separated by \cite[II.2.3 or II.2.5]{moerdijk2000proper}, hence we can apply proposition \ref{hyperconexseparated} internally in $\Lcal$ (which is boolean) and this proves that the map from $\topos$ to $\Lcal$ is atomic and internally with finite atoms. 

We need to externalise this result: Let $X$ be any object of $\topos$. Internally in $\Lcal$, the object $X$ still corresponds to an object of $\topos$ (now seen as a topos internally in $\Lcal$). Because of the first observation, it is true internally in $\Lcal$ that $X$ is a union of finite subobject (even finite atoms). Externally this mean that $X$ can be covered by subobjects $S$ defined over an open sublocale $U\subset \Lcal$ with the map $S\rightarrow U$ finite (such a $S$ corresponds to a section over $U$ of the object of finite subobject of $X$). But as $U$ is complemented, $S$ is finite in $\topos$: internally in $\topos$ if $U$ then $S$ is finite but if $\neg U$ then $S$ is empty hence also finite. This proves that any object of $X$ can be covered by finite objects and this concludes the proof.
}
}

\block{We now moves to the other half of theorem \ref{sep=locfinit}. But we need first the following lemma, and we will need it to be valid in intuitionist mathmatics.

\Lem{\label{lemmafiniteimpsep}Let $\topos$ be a (possibly non-boolean) topos which admit a generating family of of objects $X$ which are finite in the sense that $\exists n, X \simeq \{1,\dots,n \}$ holds internally. Let $p$ and $q$ be two points of $\topos$. Then the locale $I$ of isomorphisms from $p$ to $q$, defined as the pullback:

\[
\begin{tikzcd}[ampersand replacement=\&]
I \arrow{r} \arrow{d} \& * \arrow{d}{p} \\
* \arrow{r}{q} \& \topos \\
\end{tikzcd}
\]

is compact.

Moreover, this lemma holds in intuitionist mathematics.
}

\Dem{This locale $I$ classifies the theory of isomorphisms from $p$ to $q$. By hypothesis, one can construct a site $\Ccal$ for $\topos$ whose representable objects corresponds to internally finite objects. As point, $p$ and $q$ can then be seen as flat continuous functor from $\Ccal$ to sets (see \cite[VII.10]{maclane1992sheaves}), and a morphism from $p$ to $q$ can then be describe as a collection of map from $p(c)$ to $q(c)$ for $c \in \Ccal$ satisfying the natural transformation condition. Because of the assumption, all the $p(c)$ and $q(c)$ are finite sets, hence the theory of collection of map from $p(c)$ to $q(c)$ is compact by the localic Tychonof theorem and the condition of being a natural transformation can be written as a series of equality which, because equality is decidable in finite objects with this notion of finiteness, is going to be an intersection of closed sublocale hence is again closed and is hence compact.

An isomorphism is given by a couple of morphisms which are inverses of each other hence the classifying locale for isomorphisms is a closed sublocale of the locale hom$(p,q) \times$hom$(q,p)$ and hence is also compact, which concludes the proof.
}

}

\block{\label{sepLEM2}We are now ready to concludes the proof of theorem \ref{sep=locfinit}, we can even do slightly more general:

\Cor{Let $\topos$ be a topos admitting a generating family of object which are internally finite (in the sense that internally $\exists n, X \simeq \{1,\dots,n\}$) then $\topos$ is separated. Moreover this corollary holds in constructive mathematics. }

As a finite object in a boolean topos automatically satisfy this strong from of finiteness (because under the law of excluded middle all these possible definitions of finite sets are equivalent) this automatically imply the last part of the theorem

\Dem{The key observation, is that the hypothesis of being generated by a familly of finite object is pullback stable, in the sense that if $\topos$ satisfy it and $\Ecal$ is another topos then $\topos \times \Ecal$ satisfies it internally in $\Ecal$ because a pullback of a finite object is finite and a site for $\topos \times \Ecal$ internally in $\Ecal$ is given by pulling back to $\Ecal$ a site for $\topos$ in sets. Hence as the lemma \ref{lemmafiniteimpsep} has been proven constructively its conclusion will holds internally in $\Ecal$ for all the the pullback of $\topos$. The trick is then as usual to apply this to the universal case: that is internally in $\Ecal = \topos \times \topos$ in which (the pullback of) $\topos$ has two canonical point, namely $\pi_1$ and $\pi_2$ and the locale of isomorphismes between them is exactly the diagonal embedings of $\topos$ in $\topos \times \topos$, hence one get that this map is proper and hence that $\topos$ is separated. 
}

Note that even when $\topos$ is boolean, $\topos \times \topos$ is in general not (essentially unless $\topos$ is atomic) hence we really needed lemma \ref{lemmafiniteimpsep} to be proved in intuitionist mathematics in order to obtain corollary \ref{sepLEM2} even in classical mathematics.
}

\block{We conclude this section with the notion of locally separated topos.
\Def{An object $X$ of $\topos$ is said to be separating if the slice topos $\topos_{/X}$ is separated. A topos is said to be locally separated if it has a generating family of separating objects.}

A slice of a separated topos (by a decidable object) is again separated, hence in a boolean topos as soon as one has an arrow $X \rightarrow Y$ with $Y$ separating, $X$ is again separating. 
As moreover a coproduct of separating object is separating, a boolean topos (or more generally a locally decidable topos) is locally separated if an only if the terminal object an be covered by separating object.}

\block{\label{lemmafinitesepratingcover}In addition, the class of separating objects also enjoy another stability property:

\Lem{Let $f:X \twoheadrightarrow Y$ be a surjection in $\topos$. Assume that $f$ is a finite map, then $Y$ is separating.}

\Dem{The fact that $f$ is a finite surjection means that the induced map $\topos_{/X} \rightarrow \topos_{/Y}$ is a proper surjection (indeed, internally in $\topos_{/Y}$ this means that the object $X \rightarrow Y$ is inhabited and finite). The result then follows immediately from \cite[Prop II.2.1.(iii)]{moerdijk2000proper}, asserting that if $f\circ g$ is separated with $g$ a proper surjection then $f$ is separated. }

}

\section{Main result and examples}
\label{sectionmain}

\block{\label{COR_IM_imp_locsep}\Lem{Let $\topos$ be a boolean topos endowed with an invariant measure which is well supported on the terminal object, then $\topos$ is separated.}

\Dem{By (IM1), $\topos$ is generated by objects with $\mu(X)<\infty$. For such an object $X$, let $P$ be the suboject of $\un$ corresponding to the proposition ``$X$ is infinite". The map from $X \times P$ to $P$ is $\infty$-to-$1$, hence by $(IM3)$ one has $\mu(P)=0$. As $\mu$ is supposed to be well supported on $\un$ this imply that $P=0$ hence $X$ is internally finite, and hence $\topos$ is generated by finite objects and hence separated by theorem \ref{sep=locfinit}.}

Note that a local application of this result automatically show that if $X$ is an object of a topos $\topos$ endowed with an invariant measure which is well supported on $X$ then $X$ is a separating object. In particular $\topos$ is automatically locally separated. We will improve this with theorem \ref{mainth}.
}

\block{\label{invmeasonseparated}\Prop{Let $\topos$ be a boolean separated topos, and $\mu$ a well supported locally finite valuation on $\sub(\un)$, then there is an invariant measure $\tilde{\mu}$ one $\topos$ defined by:

\[ \tilde{\mu}(X) = \int_{\un} |X| d\mu,  \]

Moreover $\mu$ is the restriction of $\tilde{\mu}$ to subobjects of $\un$ and every invariant measure on $\topos$ is of the form $\tilde{\mu}$ for some locally finite well supported valuation on $\sub(\un)$.
}
\Dem{
We will first show that $\tilde{\mu}$ satisfies the three points $(IM1)-(IM3)$ of the definition \ref{invmeasDef}.

\begin{enumerate}
\item[(IM1)] One easily sees that $\tilde{\mu}$ is a valuation on $\sub(X)$ for any $X \in |\topos|$ essentially because the cardinal is internally a valuation on $X$ and that the integral of an internal valuation is again a valuation because the notion of integral we are using is linear and preserve directed supremums of positive functions.

It remains to prove that this valuation is locally finite. If $X$ is a finite object of topos, then $|X|$ is everywhere finite hence we can cover $\un$ by subobjects $U$ such that $\int_U |X|<\infty$, and for such $U$, the object $X_U=X \times U$ satisfy $\tilde{\mu}(X_U) = \int_U |X| <\infty$ and as the $X_U$ form a covering of $X$ and that every object of $\topos$ can be covered by finite objects this concludes the proof.

\item[(IM2)] Let $X$ such that $\tilde{\mu}(X)=0$, then as $\mu$ is well supported this imply that $|X|=0$ hence $X$ is the zero object. This proves that $\tilde{\mu}$ is well supported on every object.

\item[(IM3)] if $f:X\rightarrow Y$ is a $n$-to-$1$ map then internally $|X|=n|Y|$. hence if $\int |X|<\infty$, as $\mu$ is well supported, then $|X|<\infty$ everywhere, and hence one indeed get $\int |Y| = \int |X|/n$.

\end{enumerate}

Now, if $U \in \sub(\un)$ then $|U|$ is just the characteristic function of $U$ hence $\tilde{\mu}(U) = \mu(U)$.

Conversely, if $\mu$ is an invariant measure on $\topos$ and if $X$ is any object such that $\mu(X)<\infty$ then by proposition \ref{IMintegral} applied to $f:X\rightarrow \un$ (the terminal object) and $h=1$ one has:

\[ \mu(X) = \int_X 1 d\mu = \int_{\un} \left( \sum_{x \in X} 1 \right) d\mu = \int_{\un} |X| d\mu \]

Hence  $\mu$ is of the announced form. We just have to check that the measure induced by $\mu$ on $\sub(\un)$ is well supported, but if $v \subset \un$ is of zero measure any object defined over $V$ has measure zero because of the previous formula, hence axiom $(IM2)$ imply that $V$ is empty.

}}

\block{\label{mainth}\Th{Let $\topos$ be a boolean topos then $\chi$ is well supported (i.e. internally inhabited) if and only if $\topos$ is integrable and locally separated. Moreover, if $\topos$ is integrable and locally separated then $\chi$ is a principal bundle for the action of $\Rt$.}

\Dem{First let $\topos$ be a boolean topos such that $\chi$ is well supported. We will prove that $\topos$ is integrable and locally separated.

Let $X\rightarrow \chi$ be any map, then by definition $\topos_{/X}$ is endowed with an invariant measure. In particular $\topos_{/X}$ admit a generating family of object $Y$ on which the measure is well supported, and hence by lemma \ref{COR_IM_imp_locsep} a generating family of separating object of $\topos_{/X}$. Such object are also separating in $\topos_{/X}$, and as $\chi$ is assumed to be well they form a generating family of $\topos$. This proves that $\topos$ is locally separated. Moreover, as we have constructed a generating family of objects of $\topos$ which admit well supported measures this also proves that $\topos$ is integrable.

\bigskip

Conversely assume that $\topos$ is integrable and locally separated, we will prove that $\chi$ is well supported and is a principal $\Rt$-bundle.

Let $X$ be any separating object of $\topos$. As $\topos$ is integrable one can cover $X$ by subobject $S \subset X$ such that $S$ is endowed with a well supported measure. Such an object $S$ is separating and endowed with a well supported measure, hence there is an invariant measure on $\topos_{/S}$ by proposition \ref{invmeasonseparated}, and hence a map from $S$ to $\chi$. As such objects from a generating family of $\topos$, they in particular form a covering of $\un$ and hence $\chi$ is inhabited.

Moreover, still by proposition \ref{invmeasonseparated} any two functions from such an object $S$ to $\chi$ correspond to two measures on $S$, and hence (as $\sub(S)$ is boolean) the Radon-Nikodym theorem for complete boolean algebras imply that they differ only by multiplication by a function from $S$ to $\Rt$. As objects of this form are generating this exactly proves that $\chi$ is a $\Rt$-principal bundle.
}
}

\block{Let us now give some examples of boolean integrable locally separated toposes and see what is this modular bundle.

\bigskip

a)First if $\topos$ is a topos of sheaves over a boolean locale $\Lcal$ then $\topos$ is always separated and $\chi$ is the sheaf of locally finite well supported measure on $\Lcal$. It is inhbited if and only if $\Lcal$ is integrable and in this case it always has a section (assuming the axiom of choice, every inhabited sheaves over a boolean locale has a section). Moreover, it is a principal bundle exactly because of the Radon-Nikodym theorem.

\bigskip

b) If $\topos$ is the topos of $G$-sets for $G$ a discrete group. Let $G$ be endowed with its left action on itself then $\topos_{/G}$ is isomorphic to the topos of sets. Hence $G$ is separating, and $\chi$ is as a set $\R^{>0}$ (indeed, element of $\chi$ are morphism from $G$ to $\chi$, i.e. invariant measure on $\topos_{/G} \simeq $ sets. On then easily check that the $G$ action on $\chi$ is trivial, and that $\topos$ do have an invariant measure given by: 

\[ \mu(X) = \sum_{ x \in (X/G)} \frac{1}{|\text{Stab}(x)|} \]

Finally, $\topos$ is separated if and only if $G$ is finite, in which case the previous formula can be rewritten as $|X|/|G|$.

\bigskip

c) More generally, If $\topos$ is the topos of $G$ equivariant sheaves over a boolean locale $\Lcal$ endowed with an action of a discrete group $G$, then there is an object $X$ such that $\topos_{/X}$ is isomorphic to the topos of sheaves over $\Lcal$. In particular morphism from $X$ to $\chi$ are exactly measure on $\Lcal$, and a short calculation involving this object $X$ allows to deduce from this that $\chi$ is the sheaves of measure on $\Lcal$ endowed with the natural action of $G$ on it. 

In particular, an invariant measure on $\topos$ is exactly a well suported locally finite measure on $\Lcal$ invariant by the action of $G$ (this is where the name ``invariant measure" comes from). Also note that the formula giving the invariant measure on $\topos$ in terms of the $G$-invariant measure on  $\Lcal$ might be complicated in the more general case. This can be generalized easily to any boolean etendu.

\bigskip

d) Consider now the case where $G$ is a etale complete localic group\footnote{A reader unfamiliar with this terminology can replace ``etale complete localic group" by pro-discrete topological group.} and $\topos$ is the topos of continuous $G$-sets. On can then see that $\topos$ is always integrable, it is separated if and only if $G$ is compact and locally separated if and only if $G$ is locally compact. A short computation yields that $\chi$ is simply $\R^{>0}$ endowed with the action of $G$ trough its modular character. In particular, $\topos$ has an invariant measure if and only if $G$ is unimodular and in this case the invariant measure is given by the same formula as for discrete group but by replacing $|\text{Stab}(x)|$ by $\mu(\text{Stab}(x))$ for a bi-invariant\footnote{The reader can observe that as $x$ is an element of $X/G$ and not of $X$ this only make sense if $\mu$ is bi-invariant.} Haar measure $\mu$.

\bigskip

e) In the general case where $\topos$ is the topos of equivariant sheaves over an etale complete groupoid with a boolean space of objects, then one can think of $\chi$ as the sheave of right invariant left Haar systems endowed with the left action of the groupoid, and hence an invariant measure on $\topos$ can be thought of as a bi-invariant Haar system. But this cannot be made into a precise statement without developing first a theory of Haar systems well adapted to the localic framework.
 }

\section{Time evolution of Hilbert bundles and relation to the modular theory.}

\block{From the modular principal bundle $\chi$ one can define the time evolution of Hilbert spaces in the following way:

$\chi$ is $\Rt^{>0}$ bundle, and for each $t$, as $x \mapsto x^{it}$ define a continuous morphism from $\R^{>0}$ to the circle group, one can define for each $t$ a principal bundle $\chi^{it}$ for the circle group $\chi^{it}$. A principal bundle for the circle group is the same as a internal one dimensional Hilbert space, hence for each $t$ one has a internal one dimension Hilbert space $F_t$ in $\Tcal$. One can then easily see that because $\chi^{it}\otimes \chi^{it'} = \chi^{i(t+t')}$ one has $F_{t}\otimes F_{t'}=F_{t+t'}$. Hence $\sigma_t : \Hcal \mapsto \Hcal \otimes F_t$ defines a one parameter family of endofunctor of the $W^{*}$-category of $\Tcal$-Hilbert space.
}

\block{One can give a more explicit formula for these constructions:
For $t \in \R$, $F_t$ is the sheaf defined as:

\[ F_t = \{ f:\chi \rightarrow \C | \forall \alpha \in \chi, \forall r \in \R^{>0}, f(r.\alpha) = r^{-it} f(\alpha) \} \] 

The addition, multiplication and scalar product are defined pointwise, and the internal choice of a $v \in \chi$ induce an internal isomorphism\footnote{this does not mean that $F_t$ is trivial, but only that it is ``locally trivial of rank one".} between $F_t$ and $\C$ preserving all these structures.

The isomorphism between $F_t \otimes F_t'$ and $F_{t+t'}$ is given by sending $u \otimes v$ to the pointwise product of $u$ and $v$.

}

\block{The time evolution is then defined by: $\sigma_t(\Hcal) = F_t \otimes \Hcal $ for $\Hcal$ a $\topos$-Hilbert space. Equivalently, one can define:

\[ \sigma_t \Hcal =  \{ f:\chi \rightarrow \Hcal | \forall \alpha \in \chi, \forall r \in \R^{>0}, f(r.\alpha) = r^{-it} f(\alpha) \}  \]

The $\sigma_t$ are endofunctors of the category of $\topos$-Hilbert spaces.

}

\block{The last part of this paper is devoted to relate the time evolution to the modular time evolution of von Neumann algebra, and to explain how global sections of $\chi$ give rise to trace. It is actually not true that this family of functors $\sigma_t$ is the modular time evolution of the full $W^{*}$-categories of $\Tcal$-Hilbert spaces. This category is too big and there are Hilbert spaces in it which have nothing to do with the geometry of $\topos$. In order to get a result in this spirit we need to focus on a smaller, more reasonable full subcategory, essentially corresponding to the notion of ``square-integrable" representations of a group. More precisely, we will focus our attention the subcategory which is generated (in the sense of the generators of \cite[Proposition 7.3]{wstarcat}) by $\topos$-Hilbert spaces of the form $l^{2}(X)$ for $X$ a separating object, and hence it is enough to study the Hilbert spaces of the form $l^{2}(X)$ for $X$ a separating object.}

\block{For now on, $X$ denotes a \emph{separating} object of a boolean integrable locally separated topos $\topos$. We also choose $\mu$ a locally finite well supported measure on $X$, hence an invariant measure on $\topos_{/X}$, which corresponds to a morphism $\lambda$ from $X$ to $\chi$.

\bigskip

We denote by $l^2(X)$ the Hilbert space internally defined as the set of square summable $X$-indexed sequences, with (internally) its generator $(e_x) \in l^2(X)$ for each $x \in X$. We also denote by $A=\Bcal(l^2(X))$ the (external) von Neumann algebra of globally bounded operators on $l^2(X)$.
}

\block{
$\mu$ can be used to construct a normal locally finite weight on $A$, also denoted $\mu$ defined by:

\[ \forall f \in A^+,  \mu(f) = \int_{x \in X} \scal{x}{fx} d\mu \]

$\mu$ is locally finite because if $V$ is a subset of finite measure of $X$, and if $P_V \in A$ denotes the orthogonal projection on $l^2(V) \subset l^2(X)$ then $P_V f P_V$ has measure smaller than the measure of $V$ times the norm of $f$ and (assuming $f$ is positive) when $V$ vary among all finite measure subsets of $X$, this constitute a directed net of operator whose supremum if $f$.
}

\block{$\lambda$ can be used to construct an isomorphism $\phi_t: l^2(X) \simeq \sigma_t l^2(X)$.

Indeed, one can define (internally):
\[ \phi_t(e_x) := \left( \alpha \mapsto \left(\frac{\alpha}{\lambda(x)}\right)^{-it} e_x \right) \in \sigma_t l^{2}(X) \]

where $\alpha/\lambda(x)$ denotes the unique element $r$ of $\R^{>0}$ such that $\alpha = r. \lambda(x)$. Defined this way $\phi_t(e_x)$ indeed satisfies the relation $\phi_t(e_x)(r.\alpha)=r^{-it} \phi_t(e_x)$, showing that it is an element of $\sigma_t (l^2(X))$. Moreover the $\phi_t(e_x)$ are of norm one and pairwise orthogonal, hence $\phi_t$ indeed defines an isometric map $l^2(X) \rightarrow \sigma_t (l^2(X))$.

It satisfies in particular, $\phi_0 = Id$ and $\sigma_t(\phi_{t'}) \circ \phi_t = \phi_{t+t'}$ hence $\sigma_{t}\phi_{-t}$ constitutes an inverse for $\phi_t$ showing that it is an isomorphism.

}

\block{\label{timeevolutionofA}
Finally, as $l^2(X)$ is ``fixed" by the time evolution (as attested by the isomorphism $\phi_t$) one obtains an action of $\R$ directly on $A$, via:

\[\forall a\in A, \theta_t(a) = \phi_{t}^{-1} \sigma_t(a) \phi_t \]

One easily checks that it is an action of $\R$ (either directly or from the following proposition). This time evolution on $A$ can be more explicitly described on the matrix elements by:

\Prop{ For $a$ an element of $A$ and $x,y$ internal elements of $X$: 
\[ \scal{e_y}{\theta_t(a) e_x} = \left(\frac{\lambda(y)}{\lambda(x)}\right)^{-it} \scal{e_y}{a e_x} \]}

With $\frac{\lambda(y)}{\lambda(x)}$ denoting the unique element $r(x,y)$ of $\R^{>0}$ such that $r(x,y)\lambda(x)=\lambda(y)$.

\Dem{One has by definition of $\theta$ and as the $\phi_t$ are isometric:

\[\scal{e_y}{\theta_t(a) e_x} =  \scal{\phi_t(e_y)}{ \sigma_t(a) \phi_t(e_x)} \]
But:
\[ \phi_t(e_x) = \alpha \mapsto \left(\frac{\alpha}{\lambda(x)} \right)^{-it} e_x \]
And similarly for $y$, hence:
\[\begin{array}{r c l}
\scal{\phi_t(e_y)}{ \sigma_t(a) \phi_t(e_x)} &=&\displaystyle \scal{\left(\frac{\alpha}{\lambda(y)} \right)^{-it} e_y}{\left(\frac{\alpha}{\lambda(x)} \right)^{-it} a e_x}\\
& =&\displaystyle \left(\frac{\lambda(y)}{\lambda(x)}\right)^{-it} \scal{e_y}{a e_x} 
\end{array} \]

}

}

\block{\label{KMScond}Finally all these structures satisfy the modular, or Kubo-Martin-Schwinger, condition:

\Prop{For each $u \in A_+$ such that $\mu(u)$ is finite one has $\mu(\theta_t(u))=\mu(u)$.

Let $u,v \in A$ such that $\mu(u^*u)$,$\mu(uu^*)$,$\mu(v^*v)$ and $\mu(vv^*)$ are all finite. Then there exists a complex function $F_{u,v}$ defined on $\{z \in \C | Im(z) \in [-1,0] \}$ and holomorphic on its interior such that for all real numbers $t$:

\[ F_{u,v}(t)=\mu(\theta_t(u) v) \hspace{30pt} F_{u,v}(t-i)=\mu(v \theta_t(u)) \]

}

This proves that $\theta_t$ is indeed the modular group of automorphisms of the algebra $A$, associated to the semi finite normal weight $\mu$. See \cite[Chapter VIII]{takesaki2003theoryII}. Moreover, if $\topos$ is endowed with an invariant measure, then $\chi$ as a global section, $\lambda$ can be chosen to be constant equal to this global section. In this case, the formula \ref{timeevolutionofA} shows that $\theta_t$ is the identity for all $t$ and hence this last result shows that $\mu$ is a normal semi-finite trace on $A$ constructed out of the invariant measure.

\Dem{From the formula given in \ref{timeevolutionofA} one can see that $\theta_t$ left unchanged the diagonal coefficients of $u$, and as $\mu$ is defined as the integral of the diagonal coefficients one immediately has that $\mu(\theta_t(u))=\mu(u)$.

Let, for any $a\in A$, the function $a_x^y$ of matrix coefficients be defined internally by $a_x^y=\scal{e_y}{ae_x}$ (it is a function on $X \times X$). A general formal computation gives that for $a,b \in A$:

\[ \begin{array}{r c l}\mu(ab) &= &\displaystyle \int_{x \in X} \scal{e_x}{ ab e_x} d \mu \\
&=&\displaystyle \int_{x \in X} \scal{e_x}{\sum_{y \in X} a b_x^y e_y } d\mu  \\
&=&\displaystyle \int_{x \in X} \sum_{y \in X} b_x^y a_y^x d\mu \\
&=&\displaystyle \int_{(x,y) \in X^2} b_x^y a_y^x d \pi_1^* \mu 
\end{array}
\]

The last equality, corresponds to proposition \ref{IMintegral}, and holds only if $b_x^y a_y^x$ is a positive function, or if the integral is finite when we replace $b_x^y a_y^x$ by $|b_x^y a_y^x|$.

In particular, it holds when $a,b$ are $(u^*,u)$,$(u,u^*)$,$(v^*,v)$ or $(v,v^*)$, and this together with the finiteness hypothesis on $u$ and $v$ shows that all the four integrals of $|a_x^y|^2$ and $|b_x^y|^2$ with respect to both $d\pi_1^*\mu$ and $d\pi_2^*\mu$ on $X \times X$ are finite. Also note that under the correspondence between measures on $X \times X$ and functions from $X \times X$ to $\chi$, the measure $\pi_1^*\mu$ corresponds to the function $(x,y)\mapsto \lambda(x)$ and $\pi_2^*\mu$ to the function $(x,y)\mapsto \lambda(y)$. Hence one has:

\begin{equation}\label{relpi*} \pi_2^* \mu = \left(\frac{\lambda(y)}{\lambda(x)}\right) \pi_1^*\mu \end{equation}

For any complex number $z$ such that $Im(z) \in [-1,0]$ one has:

\begin{multline} \label{domintegrand} \left| \left(\frac{\lambda(y)}{\lambda(x)}\right)^{iz} u^x_y v_x^y \right|= \left(\frac{\lambda(y)}{\lambda(x)}\right)^{-Im(z)} |u^x_y||v_x^y| \\ \leqslant |u^x_y|^2 + |v_x^y|^2+\left(\frac{\lambda(y)}{\lambda(x)}\right) |u^x_y|^2+\left(\frac{\lambda(y)}{\lambda(x)}\right)|v_x^y|^2  \end{multline}

And the four terms on the right have a finite integral on $X \times X$ with respect to $\pi_1^*\mu$ (because of (\ref{relpi*}) for the last two), hence one can define the following function $F_{u,v}$ which is finite and continuous for $Im(z) \in [-1,0]$ and the previous formal computation can be applied to it.

\[ F_{u,v}(z) = \int_{(x,y) \in X \times X} \left(\frac{\lambda(y)}{\lambda(x)}\right)^{iz} u^x_y v_x^y d\pi_1^*\mu \]

Putting together the formal computation done for $\mu(ab)$ and the expression given in \ref{timeevolutionofA} for the matrix coefficients of $\theta_t(a)$, one has for $t$ real $F_{u,v}(t) = \mu(\theta_t(u)v)$, and using equation (\ref{relpi*}) one gets that:

\[\begin{array}{rcl} F_{u,v}(t-i) &= & \int_{(x,y) \in X \times X} \left(\frac{\lambda(y)}{\lambda(x)}\right)^{it} u^x_y v_x^y  \frac{\lambda(y)}{\lambda(x)}d\pi_1^*\mu \\ 
&=& 
\int_{(x,y) \in X \times X} \left(\frac{\lambda(y)}{\lambda(x)}\right)^{it} u^x_y v_x^y d\pi_2^*\mu \\
&=& 
\int_{(y,x) \in X \times X} \left(\frac{\lambda(y)}{\lambda(x)}\right)^{it} u^x_y v_x^y d\pi_1^*\mu \\
&=& \mu(v \theta_t(u))
\end{array} \]

It remains to be proven that $F_{u,v}$ is holomorphic. Let $V_n$ be the subobject of $X \times X$ on which the function $\frac{\lambda(y)}{\lambda(x)}$ is between $(1/n)$ and $n$. one has $\bigcup V_n = X \times X$ and consider:

\[ F^n_{u,v} = \int_{(x,y) \in V_n} \left(\frac{\lambda(y)}{\lambda(x)}\right)^{iz} u^x_y v_x^y d\pi_1^*\mu  \]

The functions $F^n_{u,v}$ are holomorphic, and the inequality (\ref{domintegrand}) shows that $F^n_{u,v}$ converges to $F_{u,v}$ uniformly in $z$ on all its domain of definition, showing that $F_{u,v}$ is holomorphic on the interior of its domain.

}

}

\bibliography{Biblio}{}
\bibliographystyle{plain}

\end{document}